\documentclass[11pt]{amsart}

\newtheorem{theorem}{Theorem}[section]

\theoremstyle{definition}

\theoremstyle{remark}
\newtheorem{remark}[theorem]{Remark}

\numberwithin{equation}{section}

\begin{document}

\newcommand{\spacing}[1]{\renewcommand{\baselinestretch}{#1}\large\normalsize}
\spacing{1.14}

\title{Randers Metrics of Berwald type on 4-dimensional hypercomplex Lie groups}

\author {H. R. Salimi Moghaddam}

\address{Department of Mathematics, Faculty of  Sciences, University of Isfahan, Isfahan,81746-73441-Iran.} \email{salimi.moghaddam@gmail.com and hr.salimi@sci.ui.ac.ir}


%%\date{\today}

\begin{abstract}
In the present paper we study Randers metics of Berwald type on
simply connected $4$-dimensional real Lie groups admitting
invariant hypercomplex structure. On these spaces, the Randers
metrics arising from invariant hyper-Hermitian metrics are
considered. Then we give explicit formulas for computing flag
curvature of these metrics. By this study we construct two
$4$-dimensional Berwald spaces, one of them has non-negative flag
curvature and the other one has non-positive flag curvature.
\end{abstract}

\maketitle

%%----------------------Introduction---------------------------
\section{\textbf{Introduction}}\label{Intro}
Hyper-K$\ddot{a}$hler Geometry with Torsion (HKT-geometry) is an
important field in differential geometry which have many
applications in theoretical physics. These structures appear in
supersymmetric sigma model and the study of black holes (see
\cite{GiPaSt} and \cite{Po}.). The background object of
HKT-geometry is hyper-Hermitian manifolds.\\
A special class of hyper-Hermitian manifolds are Lie groups. On
Lie groups, we can consider those hypercomplex structure and those
hyper-Hermitian (Riemannian) metrics which are invariant under
left action of $G$ on itself. The simply connected $4$-dimensional
real Lie groups admitting invariant hypercomplex structure
equipped with left invariant hyper-Hermitian (Riemannian) metrics
classified by M. L. Barberis (see \cite{Ba1} and \cite{Ba2}.).
Also we completely described the Levi-Civita connection and
sectional curvature of these metrics in \cite{Sa4}.\\
On the other hand Finsler metrics which are a generalization of
Riemannian metrics have been found many applications in physics
(\cite{AnInMa}, \cite{As}). Invariant Finsler metrics on
homogeneous spaces and Lie groups are of interesting Finsler
metrics which have been studied in the recent years (\cite{DeHo1},
\cite{DeHo2}, \cite{EsSa1}, \cite{EsSa2}, \cite{Sa1}, \cite{Sa2}
and \cite{Sa3}).\\
In this paper we study left invariant Randers metrics (a special
Finsler metric) which arise from left invariant hyper-Hermitian
(Riemannian) metrics and parallel left invariant vector fields on
simply connected $4$-dimensional real Lie groups admitting
invariant hypercomplex structure. We give the explicit formula for
computing flag curvature, a generalization of sectional curvature
for Finsler metrics, of such Randers metrics. In this way we
introduce two special complete Randers spaces of Berwald type, one
of them has non-negative flag curvature and the other one has
non-positive flag curvature.

%%---------------------Preliminaries--------------------------

\section{\textbf{Preliminaries}}\label{Prelim}

Suppose that $M$ is a $4n$-dimensional manifold. Also let $J_i,
i=1,2,3,$ be three fiberwise endomorphism of $TM$ such that
\begin{eqnarray}
  J_1J_2&=& -J_2J_1=J_3, \label{JJ} \\
  J_i^2 &=& -Id_{TM},  \ \ \ \ \ \ \ i=1,2,3, \\
  N_i &=& 0, \ \ \ \ \ \ \ i=1,2,3,
\end{eqnarray}
where $N_i$ is the Nijenhuis tensor (torsion) corresponding to
$J_i$ defined as follows:
\begin{eqnarray}
   N_i(X,Y)=[J_iX,J_iY]-[X,Y]-J_i([X,J_iY]+[J_iX,Y]),
\end{eqnarray}
for all vector fields $X,Y$ on $M$. Then the family
$\mathcal{H}=\{J_i\}_{i=1,2,3}$ is called a hypercomplex structure
on $M$.\\
In fact three complex structures $J_1, J_2$  and $J_3$ on a
$4n$-dimensional manifold $M$  form a hypercomplex structure if
they satisfy in the relation \ref{JJ} (since an almost complex
structure is a complex structure if  and only if it has no
torsion, see \cite{KoNo} page 124.).\\
A Riemannian metric $g$ on a hypercomplex manifold
$(M,\mathcal{H})$ is called hyper-Hermitian if
$g(J_iX,J_iY)=g(X,Y)$, for all vector fields $X,Y$ on $M$ and
$i=1,2,3$.\\
A hypercomplex structure $\mathcal{H}=\{J_i\}_{i=1,2,3}$ on a Lie
group $G$ is said to be left invariant if for any $a\in G$,
\begin{eqnarray}
  J_i=Tl_a\circ J_i\circ Tl_{a^{-1}},
\end{eqnarray}
where $Tl_a$ is the differential function of the left translation
$l_a$.\\

Let $M$ be a smooth $n-$dimensional manifold and $TM$ be its
tangent bundle. A Finsler metric on $M$ is a non-negative function
$F:TM\longrightarrow \Bbb{R}$ which has the following properties:
\begin{enumerate}
    \item $F$ is smooth on the slit tangent bundle
    $TM^0:=TM\setminus\{0\}$,
    \item $F(x,\lambda y)=\lambda F(x,y)$ for any $x\in M$, $y\in T_xM$ and $\lambda
    >0$,
    \item the $n\times n$ Hessian matrix $[g_{ij}(x,y)]=[\frac{1}{2}\frac{\partial^2 F^2}{\partial y^i\partial
    y^j}]$ is positive definite at every point $(x,y)\in TM^0$.
\end{enumerate}
In 1941, G. Randers introduced an important type of Finsler
metrics by using Riemannian metrics and $1$-forms on manifolds
(\cite{Ra}). Randers metrics are as follows:
\begin{eqnarray}
   F(x,y)=\sqrt{g_{ij}(x)y^iy^j}+b_i(x)y^i,
\end{eqnarray}
where $g=(g_{ij}(x))$ is a Riemannian metric and $b=(b_i(x))$ is a
nowhere zero 1-form on $M$. It has been shown that $F$ is a
Finsler metric if and only if $\|b\|=b_i(x)b^i(x)<1$, where
$b^i(x)=g^{ij}(x)b_j(x)$ and $[g^{ij}(x)]$ is the inverse matrix
of $[g_{ij}(x)]$.\\
There is another way for describing Randers metrics. We can
replace the $1$-form $b=(b_i(x))$ with its dual, therefore Randers
metrics can be defined as follows:
\begin{eqnarray}\label{F}
  F(x,y)=\sqrt{g(x)(y,y)}+g(x)(X(x),y),
\end{eqnarray}
where $X$ is a vector field on $M$ such that $\|X\|=\sqrt{g(X,X)}<1$.\\
$F$  is of Berwald type if and only if $X$ is parallel with
respect to the Levi-Civita connection induced by the Riemannian
metric $g$ (see \cite{BaChSh}.).\\
A Riemannian Metric $g$ on the Lie group $G$ is called left
invariant if
\begin{eqnarray}
  g(x)(y,z)=g(e)(T_xl_{x^{-1}}y,T_xl_{x^{-1}}z) \ \ \ \ \forall x\in
  G, \forall y,z\in T_xG,
\end{eqnarray}
where $e$ is the unit element of $G$.\\
Suppose that $g$ is a left invariant Riemannian metric on a Lie
group $G$ with Lie algebra $\frak{g}$, then the Levi-Civita
connection of $g$ is defined by the following relation
\begin{eqnarray}\label{nabla}
  2g(\nabla_UV,W)=g([U,V],W)-g([V,W],U)+g([W,U],V),
\end{eqnarray}
for any $U,V,W\in\frak{g}$, where $<,>$ is the inner product
induced by $g$ on $\frak{g}$.\\
We can define left invariant Finsler metrics similar to the
Riemannian case. A Finsler metric is called left invariant if
\begin{eqnarray}
  F(x,y)=F(e,T_xl_{x^{-1}}y).
\end{eqnarray}

The simplest way  for constructing left invariant Randers metrics
on Lie groups is the use of left invariant Riemannian metrics and
left invariant vector fields. Suppose that $G$ is a Lie group, $g$
is a left invariant Riemannian metric and $X$ is a left invariant
vector field such that $\sqrt{g(X,X)}<1$, then we can define a
left invariant Randers metric $F$ as the formula \ref{F}.\\
An important quantity which associates with a Finsler space is
flag curvature. This quantity is a natural generalization of the
concept of sectional curvature in Riemannian geometry which is
computed by the following formula:
\begin{eqnarray}\label{flag}
  K(P,Y)=\frac{g_Y(R(U,Y)Y,U)}{g_Y(Y,Y).g_Y(U,U)-g_Y^2(Y,U)},
\end{eqnarray}
where $g_Y(U,V)=\frac{1}{2}\frac{\partial^2}{\partial s\partial
t}(F^2(Y+sU+tV))|_{s=t=0}$, $P=span\{U,Y\}$,
$R(U,Y)Y=\nabla_U\nabla_YY-\nabla_Y\nabla_UY-\nabla_{[U,Y]}Y$ and
$\nabla$ is the Chern connection induced by $F$ (see \cite{BaChSh}
and \cite{Sh1}.).\\

From now we suppose that $G$ is a simply connected $4$-dimensional
real Lie group.
%%---------------------Randers Metrics of Berwald type on 4-dimensional hypercomplex Lie groups--------------------------
\section{\textbf{Randers Metrics of Berwald type on 4-dimensional hypercomplex Lie groups}}
In this paper we consider left invariant hyper-Hermitian
Riemannian metrics on left invariant hypercomplex $4$-dimensional
simply connected Lie groups. These spaces have been classified by
M. L.
Barberis as follows (for more details see \cite{Ba1}.):\\
Let $G$ be a Lie group as above with Lie algebra $\frak{g}$. She
has shown that $g$ is either Abelian or isomorphic to one of the
following Lie algebras:
\begin{enumerate}
    \item $[Y,Z]=W$, $[Z,W]=Y$, $[W,Y]=Z$, $X$ central,
    \item $[X,Z]=X$, $[Y,Z]=Y$, $[X,W]=Y$, $[Y,W]=-X$,
    \item $[X,Y]=Y$, $[X,Z]=Z$, $[X,W]=W$,
    \item $[X,Y]=Y$, $[X,Z]=\frac{1}{2}Z$, $[X,W]=\frac{1}{2}W$,
    $[Z,W]=\frac{1}{2}Y$,
\end{enumerate}
where $\{X,Y,Z,W\}$ is an orthonormal basis.\\
The case (1) is diffeomorphic to $\Bbb{R}\times S^3$ and the other
cases are diffeomorphic to $\Bbb{R}^4$ (see \cite{Ba1} and
\cite{Ba2}.).\\

Now we discuss about left invariant Randers metrics of Berwald
type which can arise from these Riamannian metrics and left
invariant vector fields on these spaces.\\

We begin with Abelian case. For this case we have the following
theorem (see \cite{Sa3}.).
\begin{theorem}
Let $G$ be an abelian Lie group equipped with a left invariant
Riemannian metric $g$ and let $\frak{g}$ be the Lie algebra of
$G$. Suppose that $X\in\frak{g}$ is a left invariant vector field
with $\sqrt{g(X,X)}<1$. Then the Randers metric $F$ defined by the
formula \ref{F} is a flat geodesically complete locally
Minkowskian metric on $G$.
\end{theorem}

Now we continue with the other four cases.\\
\textbf{Case 1.} We can compute the Levi-Civita connection by
using formula \ref{nabla} (Also you can see \cite{Sa4}.). The
Levi-Civita connection is as follows:
\begin{eqnarray}
  &&\nabla_XX = 0 \ \ , \ \ \nabla_XY = 0, \ \ \nabla_XZ = 0, \ \ \nabla_XW = 0, \nonumber\\
  &&\nabla_YX = 0 \ \ , \ \ \nabla_YY = 0, \ \ \nabla_YZ = \frac{1}{2}W, \ \ \nabla_YW = -\frac{1}{2}Z,\\
  &&\nabla_ZX = 0 \ \ , \ \ \nabla_ZY = -\frac{1}{2}W, \ \ \nabla_ZZ = 0, \ \ \nabla_ZW = \frac{1}{2}Y,\nonumber\\
  &&\nabla_WX = 0 \ \ , \ \ \nabla_WY = \frac{1}{2}Z \ \ , \ \ \nabla_WZ = -\frac{1}{2}Y, \ \ \nabla_WW = 0\nonumber.
\end{eqnarray}
A simple computation shows that the only family of vector fields
which is parallel with respect to this connection is of the form
$Q=qX$ for any $q\in\Bbb{R}$. Now let $0<\|Q\|<1$ or equivalently
let $0<|q|<1$, therefore by using these left invariant vector
fields $Q$ and formula \ref{F}, $G$ admits a family of Randers
metrics of Berwald type.\\
\textbf{Case 2.} The formula \ref{nabla} shows that the
Levi-Civita connection of the Riemannian metric of $G$ is of the
form (\cite{Sa4}):
\begin{eqnarray}
  &&\nabla_XX = -Z \ \ , \ \ \nabla_XY = 0, \ \ \nabla_XZ = X, \ \ \nabla_XW = 0, \nonumber\\
  &&\nabla_YX = 0 \ \ , \ \ \nabla_YY = -Z, \ \ \nabla_YZ = Y, \ \ \nabla_YW = 0,\\
  &&\nabla_ZX = 0 \ \ , \ \ \nabla_ZY = 0, \ \ \nabla_ZZ = 0, \ \ \nabla_ZW = 0,\nonumber\\
  &&\nabla_WX = -Y \ \ , \ \ \nabla_WY = X \ \ , \ \ \nabla_WZ = 0, \ \ \nabla_WW = 0\nonumber.
\end{eqnarray}
Similar to case 1 by a simple computation we can show that the
only parallel vector fields with respect to $\nabla$ are of the
form $Q=qW, q\in\Bbb{R}$. For constructing a non-Riemannian family
of Randers metrics of Berwald type defined by formula \ref{F}, it
is sufficient to let $0<\|Q\|<1$ or equivalently $0<|q|<1$.\\
\textbf{Case 3, and case4.} The Levi-Civita connections of case 3
and 4 are of the following forms respectively:
\begin{eqnarray}
  &&\nabla_XX = 0 \ \ , \ \ \nabla_XY = 0, \ \ \nabla_XZ = 0, \ \ \nabla_XW = 0, \nonumber\\
  &&\nabla_YX = -Y \ \ , \ \ \nabla_YY = X, \ \ \nabla_YZ = 0, \ \ \nabla_YW = 0,\\
  &&\nabla_ZX = -Z \ \ , \ \ \nabla_ZY = 0, \ \ \nabla_ZZ = X, \ \ \nabla_ZW = 0,\nonumber\\
  &&\nabla_WX = -W \ \ , \ \ \nabla_WY = 0 \ \ , \ \ \nabla_WZ = 0, \ \ \nabla_WW =
  X\nonumber,
\end{eqnarray}
and
\begin{eqnarray}
  &&\nabla_XX = 0 \ \ , \ \ \nabla_XY = 0, \ \ \nabla_XZ = 0, \ \ \nabla_XW = 0, \nonumber\\
  &&\nabla_YX = -Y \ \ , \ \ \nabla_YY = X, \ \ \nabla_YZ = -\frac{1}{4}W, \ \ \nabla_YW = \frac{1}{4}Z,\\
  &&\nabla_ZX = -\frac{1}{2}Z \ \ , \ \ \nabla_ZY = -\frac{1}{4}W, \ \ \nabla_ZZ = \frac{1}{2}X, \ \ \nabla_ZW = \frac{1}{4}Y,\nonumber\\
  &&\nabla_WX = -\frac{1}{2}W \ \ , \ \ \nabla_WY = \frac{1}{4}Z \ \ , \ \ \nabla_WZ = -\frac{1}{4}Y, \ \ \nabla_WW =
  \frac{1}{2}X\nonumber.
\end{eqnarray}
It is easy to show that these connections do not admit any
parallel left invariant vector field, therefore there is not any
left invariant Randers metric of Berwald type arising from left
invariant vector fields (by using formula \ref{F}) on these
Riemannian Lie groups.
%%---------------------Flag curvature--------------------------
\section{\textbf{Flag curvature}}
In this section we discuss about the flag curvature of invariant
Randers metrics of cases 1 and 2.\\
\textbf{Case 1.} By using Levi-Civita connection for curvature
tensor we have (\cite{Sa4}):
\begin{eqnarray}
  &&R(Y,Z)Y = -R(Z,W)W = -\frac{1}{4}Z,\nonumber\\
  &&R(Y,W)W = R(Y,Z)Z = \frac{1}{4}Y,\\
  &&R(Z,W)Z = R(Y,W)Y = -\frac{1}{4}W\nonumber,
\end{eqnarray}
and in other cases $R=0$. Now let $U=aX+bY+cZ+dW$ and
$V=\tilde{a}X+\tilde{b}Y+\tilde{c}Z+\tilde{d}W$ be two arbitrary
vectors in $\frak{g}$ then we have:
\begin{eqnarray}
  R(V,U)U=-\frac{1}{4}\{(b\tilde{c}-c\tilde{b})(cY-bZ)+(b\tilde{d}-d\tilde{b})(dY-bW)+(c\tilde{d}-d\tilde{c})(dZ-cW)\}.
\end{eqnarray}
Since $F$ is of Berwald type therefore the curvature tensor of $F$
and $g$ coincide. Suppose that $\{U,V\}$ is an orthonormal basis
for $P=span\{U,V\}$ with respect to the inner product $<,>$
induced by $g$. Now by using formula
$g_U(V_1,V_2)=\frac{1}{2}\frac{\partial^2}{\partial s\partial
t}(F^2(U+sV_1+tV_2))|_{s=t=0}$ of $F$ (for an explicit formula you
can see \cite{EsSa1}.) we have:
\begin{eqnarray}
  g_U(R(V,U)U,V) &=& \frac{1}{4}(1+aq)\{(b\tilde{c}-c\tilde{b})^2+(b\tilde{d}-d\tilde{b})^2+(c\tilde{d}-d\tilde{c})^2\}\\
  g_U(U,U) &=& (1+aq)^2 \\
  g_U(V,V)&=& 1+aq+(\tilde{a}q)^2 \\
  g_U(U,V) &=& \tilde{a}q(1+aq).
\end{eqnarray}
Now by using \ref{flag} we have:
\begin{eqnarray}
  K(P,U) &=&
  \frac{(b\tilde{c}-c\tilde{b})^2+(b\tilde{d}-d\tilde{b})^2+(c\tilde{d}-d\tilde{c})^2}{4(1+aq)^2}\geq
  0.
\end{eqnarray}
Therefore in the case 1 $(G,F)$ is of non-negative flag
curvature.\\
\textbf{Case 2.} The curvature tensor of Riemannian metric
(Finsler metric) of this case is of the form:
\begin{eqnarray}
  &&R(X,Y)X = -R(Y,Z)Z = Y,\nonumber\\
  &&R(X,Y)Y = R(X,Z)Z = -X,\\
  &&R(X,Z)X = R(Y,Z)Y = Z\nonumber.
\end{eqnarray}
In this case for any $U$ and $V$ we have:
\begin{eqnarray}
  R(V,U)U=-\{(a\tilde{b}-b\tilde{a})(aY-bX)+(a\tilde{c}-c\tilde{a})(aZ-cX)+(b\tilde{c}-c\tilde{b})(bZ-cY)\},
\end{eqnarray}
Let $P=\{U,V\}$ be as case 1. Therefore for the Randers metric $F$
described in case 2 we have:
\begin{eqnarray}
  g_U(R(V,U)U,V) &=& -(1+dq)\{(a\tilde{b}-b\tilde{a})^2+(a\tilde{c}-c\tilde{a})^2+(b\tilde{c}-c\tilde{b})^2\}\\
  g_U(U,U) &=& (1+dq)^2 \\
  g_U(V,V)&=& 1+dq+(\tilde{d}q)^2 \\
  g_U(U,V) &=& \tilde{d}q(1+dq).
\end{eqnarray}
Hence for the flag curvature we have:
\begin{eqnarray}
  K(P,U) &=&
  \frac{-\{(a\tilde{b}-b\tilde{a})^2+(a\tilde{c}-c\tilde{a})^2+(b\tilde{c}-c\tilde{b})^2\}}{(1+dq)^2}\leq
  0,
\end{eqnarray}
which shows that $(G,F)$ is of non-positive flag curvature.
\begin{remark}
The metrics constructed in case 1 and 2 are complete and
$J_i$-invariant, $i=1,2,3$.
\end{remark}
\begin{proof}
These metrics clearly are $J_i$-invariant, $i=1,2,3$ so we prove
the completeness. Since the metric $F$ (in cases 1 or 2) is of
Berwald type therefore the geodesics of $F$ and $g$ coincide. On
the other hand $(G,g)$ is a homogeneous Riemannian manifold, hence
$(G,g)$ is geodesically complete (see \cite{BeEhEa} page 185.).
Therefore $(G,F)$ is geodesically complete. Now the Hofp-Rinow
theorem and connectedness of $G$ will complete the proof.
\end{proof}

%%-------------------- BIBLIOGRAPHY------------------------

\bibliographystyle{amsplain}

\end{document}